\documentclass[12pt]{article}
\usepackage{amsmath,mathrsfs,amssymb}
\usepackage{amscd}
\begin{document}
\title{\bf \large REDUCTION OF QUASI-LATTICES TO LATTICES}
\author{{\sc \footnotesize C. Ganesa Moorthy $^a$} {\footnotesize and} {\sc \footnotesize  SG. Karpagavalli $^{b,}$ \footnote {Corresponding author}}}
\date{}
\maketitle
\begin{center}
{\footnotesize \it  $^a$ Department of Mathematics, Alagappa University, Karaikudi-630 004,\\ Tamil Nadu, India.\\}
{\footnotesize \it $^b$ Department of Mathematics, Vidhyaa Giri College of Arts and Science, Puduvayal-630 108, Tamil Nadu, India.}
\end{center}
\begin{center}
{{\footnotesize E-mail:}  {\tt \footnotesize ganesamoorthyc@gmail.com} (\footnotesize C. Ganesa Moorthy)\\ {\tt \quad sgkarpa@gmail.com} (SG.Karpagavalli)}
\end{center}
\newtheorem{defn}{Definition}[section]
\newtheorem{lem}[defn]{Lemma}
\newtheorem{thm}[defn]{Theorem}
\newtheorem{cor}[defn]{Corollary}
\newtheorem{rem}[defn]{Remark}
\newtheorem{exam}[defn]{Example}
\newtheorem{pros}[defn]{Proposition}
\newtheorem{note}[defn]{Notation}
\begin{abstract}
	Quasi-lattices are introduced in terms of $`$join' and $`$meet' operations. It is observed that quasi-lattices become lattices when these operations are associative and when these operations satisfy $`$modularity' conditions. A fundamental theorem of homomorphism proved in this article states that a quasi-lattice can be mapped onto a lattice when some conditions are satisfied.
\end{abstract}
{\bf Key words:} Minimal upper bound, Congruence relation, Partition.\\
{\bf AMS Subject Classification(2010):} 06B10,18B35,06C05.  
\section{Introduction}
 	The concept of a minimal upper bound is not widely known. A lattice is a partially ordered set (poset) in which any two elements have a least upper bound and a greatest lower bound. A quasi-lattice is a poset in which any two elements have a minimal upper bound and a maximal lower bound. Every quasi-lattice is a lattice. This article tries to establish fundamental facts about quasi-lattices. But, it finds that associativity of $`$meet' and $`$join' operations of quasi-lattices is a unique property of lattices. Similarly it is established that $`$modularity' is also a unique property of lattices. A fundamental theorem of homomorphism found in this article also reduces quasi-lattices into lattices. The books [3] and [2]  are referred to fundamental definitions and properties for posets and lattices. Although there are many recent articles (see, for example [4, 5, 6]) the results of these articles will not be extended to quasi-lattices, because quasi-lattices reduce to lattices when some fundamental properties are assumed.\\
\indent	A partial order  $\leq$   on a non empty set $P$ is a relation that is reflexive, anti-symmetric and transitive. A poset  $(P,\leq)$ is a non empty set $P$ with a partial order $\leq$. An element $a$ in a partially ordered set  $(P,\leq)$ is a maximal lower bound of a non empty subset $A$ of  $P$ if $a \leq x$, $\forall x \in  A$, and if there is no element $d$ in $P$ such that $a <d \leq x$, $\forall x \in  A$.Dually a minimal upper bound is defined. A partially ordered set $(P,\leq)$ is called quasi-lattice, if any two elements of $P$ have a minimal upper bound and a maximal lower bound. However, two elements in a quasi-lattice may have more than one maximal lower bound and may have more than one minimal upper bound. Let us use the notations $x \wedge y$ and $x \vee y$ to denote some (particular) maximal lower bound and some minimal upper bound of $x$ and $y$, respectively,   in a partially ordered set.
\begin{exam}
The Hasse diagram given in the Figure 1 represents a quasi-lattice. In this diagram the point $x \vee(y \vee z)$ represents another minimal upper bound of \{x,y\} in addition to
$x \vee y$. So, it is not a lattice. In this quasi-lattice, $(x\vee y)\vee z$ can never take $``$ the form"  $x \vee (y\vee z)$. So, associativity fails to be true.
\begin{center}
\unitlength 1mm 
\linethickness{0.4pt}
\ifx\plotpoint\undefined\newsavebox{\plotpoint}\fi 
\begin{picture}(51.25,48.5)
\put(23.25,43.5){\circle*{1.5}}
\put(13.5,31.25){\circle*{2.062}}
\put(34,31.5){\circle*{1.803}}
\put(24,19.5){\circle*{2.236}}
\put(5.25,19.75){\circle*{2.55}}
\put(46.5,20.25){\circle*{2.236}}
\put(25.5,3.5){\circle*{2}}
\multiput(22.5,42.75)(-.0336363636,-.0381818182){275}{\line(0,-1){.0381818182}}
\multiput(23.5,43)(.0336538462,-.0368589744){312}{\line(0,-1){.0368589744}}
\multiput(33.5,30.5)(-.0336538462,-.0384615385){260}{\line(0,-1){.0384615385}}
\multiput(24.75,20.5)(-.03125,-.0625){8}{\line(0,-1){.0625}}
\put(14,30.5){\line(0,1){0}}
\multiput(12.5,30.5)(-.03372093,-.046511628){215}{\line(0,-1){.046511628}}
\multiput(5.25,18.75)(.0442477876,-.0337389381){452}{\line(1,0){.0442477876}}
\multiput(26.25,3.75)(.0410526316,.0336842105){475}{\line(1,0){.0410526316}}
\multiput(34.25,31.25)(.0389221557,-.0336826347){334}{\line(1,0){.0389221557}}
\put(47.25,20){\line(0,1){0}}
\multiput(24,19)(.03333333,-.36666667){45}{\line(0,-1){.36666667}}
\put(27.75,38){\circle*{1.803}}
\multiput(6,20.25)(.0409441233,.0337186898){519}{\line(1,0){.0409441233}}
\multiput(14,30.75)(.0337171053,-.0361842105){304}{\line(0,-1){.0361842105}}
\put(29,18.75){\makebox(0,0)[cc]{y}}
\put(51.25,20.5){\makebox(0,0)[cc]{z}}
\put(2,17.25){\makebox(0,0)[cc]{x}}
\put(29.75,1.75){\makebox(0,0)[cc]{0}}
\put(22.75,49.0){\makebox(0,0)[cc]{$( x \vee y) \vee z$}}
\put(37,41){\makebox(0,0) [cc] {  $x \vee (y \vee z)$}}
\put(41.25,32.5){\makebox(0,0)[cc]{$ y \vee z$}}
\put(4.75,31.75){\makebox(0,0)[cc]{$x \vee y$}}
\put(46.75,4.5){Figure 1}
\end{picture}
\end{center}
\end{exam}
\section{Associative quasi-lattices}
	It would be difficult to derive many results for quasi-lattices, when associativity is not assumed.
\begin{defn} A quasi-lattice $(P, \leq)$ is called an associative lattice, if 
\begin{enumerate}
\item [(i)] $a\vee (b\vee c) = (a\vee b)\vee c$, and
\item[(ii)] $a\wedge(b\wedge c) = (a\wedge b)\wedge c$ hold for every $a,b,c \in P$.
\end{enumerate} 
\end{defn}                   
                 Here (i) means that if $a_1$ is a minimal upper bound of \{$b,c$\} and if $a_2$ is a minimal upper bound of \{$a,a_1$\}, then there is a minimal upper bound $a_3$ of \{$a,b$\} such that $a_2$ is a minimal upper bound of\{$a_3,c$\} and similarly; if $b_1$ is a minimal upper bound of \{$a,b$\} and $b_2$ is minimal upper bound of \{$b_1,c$\}, then there is a minimal upper bound $b_3$ of \{$b,c$\} such that $b_2$ is a minimal upper bound of \{$a,b_3$\}. This interpretation clarifies the meanings for the present notation. When this is followed, the meaning of the following proposition is unambigious.
\begin{pros} The following identities are true in a quasi-lattice $(P,\leq)$.
(A1): $a\vee a = a$; (A2): $a\wedge a = a$; (A3): $a\vee b = b\vee a$; (A4):$a\wedge b= b\wedge a$; (A5): $a\vee(a\wedge b) = a= (a\wedge b)\vee a$ ;
(A6):$a\wedge (a\vee b) = a = (a\vee b)\wedge a$; $\forall a,b \in P$.
\end{pros}
{\bf Proof:} Let us verify $a\vee (a\wedge b) =a$.
Let $a_1$ be a maximal lower bound of \{$a,b$\}, and $a_2$ be a minimal upper bound of \{$a_1,a$\}. Then $a_2 = a$ because $a_1 \leq a$. Other relations can also be verified in this way.\\
\indent  If $( P,\leq )$ is an associative quasi-lattice, then it further has the properties:
( A7) : $a\vee  (b\vee c) = (a\vee b)\vee c$ and
( A8) : $a\wedge (b\wedge c) = (a\wedge b)\wedge c$;  $\forall  a, b, c \in P$.
        It is known that the relations ( A1) to (A8) characterize a lattice, when $a\vee b$ and $a\wedge b$ are unique elements ( see Theorem 1 in Section 1 in Chapter 1 in [2] ). It is to be proved that an associative quasi-lattice should be a lattice.  For this purpose, let us introduce some changes in applications of the notations $\vee$  and $\wedge$. For a given poset $(P,\leq),A\subseteq P$ and $B\subseteq P$, let $A\vee B$ (respectively, $A\wedge B$) denote the collection of all elements of the form $a\vee b$ (respectively, 
$a\wedge b$) with $a\in A$ and  $b\in B$. So, for example, the relation $a\wedge(a\vee b) = a = (a\vee b)\wedge a$ will mean $\{a\}\wedge (\{a\}\vee \{b\}) = \{a\}
 = (\{a\}\vee \{b\})\wedge \{a\}$. Thus a poset $(P,\leq)$ is a quasi-lattice if and only if $\{a\}\vee \{b\}$ and $\{a\}\wedge \{b\}$  are non empty subsets of $P$, for any $a,b \in P$. It is a lattice if and only if $\{a\}\vee \{b\}$ and $\{a\}\wedge \{b\}$ are singleton subsets of $P$, for any $a,b \in P$.
\begin{thm}	A quasi-lattice $(P,\leq)$ is  associative if and only if it is a lattice.
\end{thm}
{\bf Proof:} Suppose $(P,\leq)$ is an associative quasi-lattice.
 Let $x,y\in P$ and $a,b \in  \{x\}\wedge \{y\}$.Then $a\leq y$, $a \leq x$ and 
$(\{x\}\wedge \{y\})\wedge \{a\} = \{x\}\wedge (\{y\}\wedge\ \{a\})= \{x\}\wedge \{a\}= \{a\}$, when $(\{x\}\wedge \{y\})\wedge\{a\}\supseteq \{a,b\}\wedge \{a\}= \{a\}\cup (\{b\}\wedge\{a\})$.Thus $\{a\}\wedge \{b\} = \{a\}$ so that $a \leq b$.
Similarly $b\leq a$ so that $a = b$.Thus $\{x\}\wedge \{y\}$ contains a unique element. Dually, $\{x\}\vee \{y\}$ contains a unique element. This proves that $(P,\leq)$ is a lattice.
\section{Modular quasi-lattices}
\begin{defn}	A quasi-lattice $(P,\leq)$ is said to be modular if
$\{x\}\vee (\{y\}\wedge \{z\}) = (\{x\}\vee\{y\})\wedge \{z\}$ whenever $x,y,z \in P$ and 
$x \leq z$.
\end{defn}
\begin{thm} A modular quasi-lattice $(P,\leq)$ is a lattice.
\end{thm}
{\bf Proof:} Fix $x,y$ in the given modular lattice $(P,\leq)$.
Let $a,b\in\{x\}\wedge\{y\}$. Then $a\leq x, a \leq y, b\leq x$, and $b\leq y$.
So, $(\{x\}\wedge \{y\})\vee\{a\} = \{a\}\vee(\{y\}\wedge\{x\})= (\{a\}\vee \{y\})\wedge \{x\}= \{y\}\wedge \{x\} = \{x\}\wedge\{y\}$, when $\{a,b\}\subseteq\{x\}\wedge \{y\}$.
So $\{a,b\}\vee \{a\} \subseteq\{x\}\wedge\{y\}$ and
hence $\{a,a\vee b\}\subseteq\{x\}\wedge \{y\}$.Thus $a\vee b\in \{x\}\wedge\{y\}$, when
$a\vee b \geq a, a\vee b\geq b$ , $a \in\{x\}\wedge \{y\}$ and $b\in\{x\}\wedge\{y\}$.
So,the  maximality of $a$ and $b$  implies that $a = a\vee b = b$. In particular, $\{x\}\wedge \{y\}$ contains atmost one point. Dually $\{x\}\vee\{y\}$ contains atmost one point. This proves the theorem.\\
\indent Associative quasi-lattices are lattices and modular quasi-lattices are lattices. So it is difficult to derive new results for quasi-lattices, because quasi-lattices with additional fundamental properties become lattices. However, one can derive fundamental results for ideals.
\begin{defn} A  subset $I$  $(\mathscr{F}, respectively)$ of  a quasi-lattice $(P,\leq)$ is called an ideal (a filter, respectively), if
\begin{enumerate}
\item[(i)] $a,b\in I \Rightarrow \{a\}\vee \{b\} \subseteq I$
\item[((i)] $a,b\in \mathscr{F}\Rightarrow \{a\}\wedge\{b\} \subseteq \mathscr{F}$, respectively) and
\item[(ii)] $a\in I,b\in P, b\leq a \Rightarrow b \in I$
\item[((ii)] $a\in \mathscr{F}, b\in P, b\geq a \Rightarrow b\in \mathscr{F}$, respectively).
\end{enumerate}
\end{defn}
\indent   An arbitrary intersection of ideals (filters) in a quasi-lattice is an ideal (a filter). The intersection of a filter and an ideal is sub quasi-lattice. Here, a sub quasi-lattice $(Q,\leq)$ of a quasi-lattice $(P,\leq)$ means that $\{x\}\vee\{y\}\subseteq Q$, and $\{x\}\wedge\{y\}\subseteq Q$, whenever $x,y \in Q$. The intersection of a filter with an ideal is a convex subset in view of the following (usual) definition.

\begin{defn} A subset $C$ of a quasi-lattice $(P,\leq )$ is said to be convex, 
if $a\in C$, whenever $x,y\in C$, $a\in P$  and $x\leq a\leq y$.
\end{defn}
\begin{note} To each  $A\subseteq P$, a quasi-lattice, let $(A]$ and $[A)$ denote the smallest ideal and the smallest filter, respectively, containing $A$. They exist in view of the previous remark.
\end{note}
\begin{pros} Let  $(P,\leq)$ be a quasi-lattice. Let $I(P)$ (respectively, $F(P)$) be the collection of all ideals (respectively, filters) of $(P,\leq)$. Then $I(P)$ (respectively, $F(P)$) is a complete lattice under the inclusion relation (respectively, inverse inclusion relation).
\end{pros}
{\bf Proof:} Let $(I_\lambda)_{\lambda\in A}$ be a collection of ideals in $P$, Then $\cap\{I_\lambda:\lambda \in A\}$ and $(\cup\{I_\lambda:\lambda \in A\}]$ are ideals  which are the greatest lower bound and the least upper bound of the given collection. A similar argument is applicable for filters.
\section{Congruence relations}
Ideals are associated with inverse image of a least element for a lattice homomorphism. A lattice homomorphism is associated with a congruence. Let us first define a congruence relation for a quasi-lattice.
\begin{defn}Let $(P,\leq)$ be a quasi-order lattice. An equivalence relation $\theta$ on $P$ is denoted by $x\equiv y$(mod $\theta$) when $x$ and $y$ are related in $P$ by $\theta$. Moreover, for subsets $A$,$B$ of $P$, the identity  $A\equiv B$ (mod $\theta$) means the following:
\begin{enumerate}
\item[(i)] to each $a\in A$, there is a   $b \in B$ such that $a\equiv b$ (mod $\theta$), and
\item[(ii)] to each $b \in B$, there is an   $a\in B$ such that $a\equiv b$ (mod $\theta$).\end{enumerate} The equivalence relation $\theta$ on $P$ is called a congruence relation, if $\{x_1\}\wedge\{y_1\}\equiv\{x_2\}\wedge\{y_2\}$(mod $\theta$), and $\{x_1\}\vee\{y_1\}\equiv\{x_2\}\vee\{y_2\}$(mod $\theta$), whenever $x_1\equiv x_2$(mod $\theta$) and $y_1\equiv y_2$(mod $\theta$) in $P$, and if $\{x\}\wedge\{y\}\subseteq[z]$,   when $z\in \{x\} \wedge \{y\}$ and $\{x\}\vee\{y\}\subseteq[z]$, when $z\in \{x\} \vee \{y\}$, for x,y,z in P, when $[z]$ refers to the equivalence class containing $z$, determined by $\theta$.\\
\indent	It is known that the collection of all partitions is a complete lattice under the $``$refinement" relation. The collection of all congruences on a lattice is a (complete) sublattice of the lattice of all partitions. In the same way(see the proof of theorem 3.9 in [1]), one can verify that the collection of all congruences on a quasi-lattice is a complete lattice and a sublattice of the lattice of all partitions.
\end{defn}
\begin{lem}	Let $(P,\leq)$ be a quasi-lattice,  and $\theta$ be a congruence relation on $P$. If $u\equiv v$ (mod $\theta$), $a\in\{u\}\wedge\{v\}$, $b\in\{u\}\vee\{v\}$, and if
$a\leq x\leq b$ , then $u\equiv x$ (mod $\theta$).\end{lem}
{\bf Proof:} Under the assumptions, we have $\{x\}=\{x\}\vee\{a\}\equiv\{x\}\vee(\{u\}\wedge\{v\})\equiv\{x\}\vee(\{u\}\wedge\{u\})                                                                 
\equiv(\{x\}\vee\{u\})$(mod $\theta$). Dually, we have $\{x\} = \{x\}\wedge\{b\}\equiv\{x\}\wedge(\{u\}\vee\{v\})\equiv\{x\}\wedge(\{u\}\vee\{u\})
\equiv\{x\}\wedge\{u\}$(mod $\theta$). So, we have $\{u\}= \{u\}\wedge(\{u\}\vee\{x\})=\{u\}\wedge\{x\}\equiv\{x\}$(mod $\theta$). This proves the lemma.
\begin{defn} Let $T:P_1\to P_2$ be a mapping from a quasi-lattice $P_1$ into a quasi-lattice $P_2$. It is said to be a q-lattice homomorphism, if $T(\{x\}\vee\{y\})=\{T(x)\}\vee\{T(y)\}$ and $T(\{x\}\wedge\{y\})=\{T(x)\}\wedge\{T(y)\},\forall x,y\in P$
\end{defn}
\begin{defn} Let $\theta$ be an equivalence relation on a quasi-lattice $(P,\leq)$. Let $[x]$ denote the equivalence class containing $x$. Let us say that $\theta$ satisfies the condition (*) if the following are true in $P$ :
\begin{enumerate}
\item [(i)]  If $[x]\neq[y]$, $x\leq z$ and $y\leq z$, then there are elements $a\in[x]$ and $b\in[y]$ and there is an element $d\in\{a\}\vee\{b\}$ such that $d\leq z$.
\item[(ii)]  If $[x]\neq[y]$, $x\geq z$ and $y\geq z$, then there are elements $a\in[x]$ and $b\in[y]$ and there is an element $d\in\{a\}\wedge\{b\}$ such that $d\geq z$.
\end{enumerate}
\end{defn}
	Let us now state  a fundamental theorem of homomorphism.
\begin{thm}	Let $(P,\leq)$ be a quasi-lattice. Let $\theta$ be a congruence relation on $P$ that satisfies (*) of definition 4.4. Let $P/\theta$  be the collection of all equivalence classes. Let $[x]$ denote the equivalence class containing $x$. Then $P/\theta$  is a lattice in which we have $[x]\wedge[y] =[x\wedge y]$ and $[x]\vee[y] =[x\vee y]$, for any elements $x\wedge y$ and $x\vee y$ in $\{x\}\wedge\{y\}$ and $\{x\}\vee\{y\}$, respectively. Also, the quotient mapping  $\pi: P\to P/\theta$ defined by $\pi(x) =[x]$, $x\in P$, is a surjective q-lattice homomorphism. On the other hand if $T: P\to L$ is a surjective q-lattice homomorphism from a quasi-lattice $P$ onto a lattice $L$, then $\{T^{-1}(a):a\in L\}$ defines a partition that leads to a congruence relation satisfying (*) of definition 4.4.\end{thm}
{\bf Proof:}\\
{\bf First Part:} Define $[x]\leq[y]$ if and only if $a\leq b$ for some $a\in [x]$ and some $b\in[y]$. Suppose $a_1\in[x]$ and $b_1\in[y]$ such that $a_1\leq b_1$. If $a_2\in[x]$, then $a_1\equiv a_2$ (mod $\theta$), $a_2\leq b_1\vee a_2$	(for any element of this type) and $\{b_1\}\vee\{a_2\}\equiv\{b_1\}\vee\{a_1\}\equiv\{b_1\}$(mod $\theta$).
If $b_2\in [y]$, then $b_1\equiv b_2$(mod $\theta$), $a_1\wedge b_2\leq b_2$, and
$\{a_1\}\wedge\{b_2\}\equiv\{a_1\}\wedge\{b_1\}\equiv\{a_1\}$(mod $\theta$). Thus, if $[x]\leq[y]$, then for any $a_1\in[x]$, there is a $b_1\in[y]$ such that $a_1\leq b_1$ and for any $b_2\in[y]$ there is an $a_2\in[x]$ such that $a_2\leq b_2$. Now let us verify that this relation in $P/\theta$ is a partial order relation. Since $x\leq x$, we have 
$[x]\leq [x], \forall x\in P$. To prove anti-symmetricity, assume that $[x]\leq[y]$ and $[y]\leq[x]$ for two elements $x,y\in P$. Then there is an element $y_1\in[y]$ such that $x\leq y_1$; and there is an element $x_1\in[x]$ such that $y_1\leq x_1$. Thus $x\leq y_1\leq x_1$ and $x_1\equiv x$ (mod $\theta$). By the previous lemma 4.2 it is concluded that $y_1\equiv x$ (mod $\theta$). This proves that $\leq$ is anti-symmetric in $P/\theta$.
To prove transitivity, assume that $[x]\leq[y]$ and $[y]\leq[z]$ for some $x,y,z\in P$.
Then there is an element $y_1\in[y]$ and there is an element $z_1\in [z]$ satisfying $x\leq   y_1\leq z_1$ so that $x\leq z_1$, So $(P/\theta ,\leq)$ is a poset.
To prove that $P/\theta$ is a lattice, consider an element $a\in\{x\}\wedge\{y\}$, for some fixed elements $x,y$. Then $a\leq x$ and $a\leq y$. So $[a]\leq [x]$ and $[a]\leq [y]$. Suppose $[b]\leq [x]$ and $[b]\leq [y]$ for some element $b$ of $P$, and assume that $[a]\leq[b]$. Then there is an element $b_1\in[b]$ such that $a\leq b_1$. There are elements $c_1\in[x]$ and $c_2\in[y]$ such that $b_1\leq c_1$ and $b_1\leq c_2$. By the condition (*) satisfied, there are elements $a_1\in[c_1]$ and $a_2\in[c_2]$ and  there is an element $d\in\{a_1\}\wedge\{a_2\}$ such that $b_1\leq d$. Since $\{a\} \equiv\{x\}\wedge\{y\}\equiv\{c_1\}\wedge\{c_2\}\equiv\{a_1\}\wedge\{a_2\}\equiv\{d\}$(mod $\theta$), we have the relation $[b]\leq[a]$. Thus $[a] = [b]$. This proves that $[x]\wedge[y]=[x\wedge y]$ for any element $x\wedge y,\forall x,y \in P$. Dually, one can prove that $[x]\vee[y] = [x\vee y]$, for any element $x\vee y,\forall x,y\in P$. So,  $P/\theta$  is a lattice. Other sub divisions of the first part are trivial.\\
{\bf Second Part:} Let $\theta$ be the equivalence relation induced by the partition $\{T^{-1}(a): a\in L\}$. The condition (*) of definition 4.4 has to be checked to complete the proof as the other sub divisions are trivial. If $T^{-1}(a)\neq T^{-1}(b)$, $a_1\in T^{-1}(a), b_1\in T^{-1}(b),d_1\geq a_1, d_1\geq b_1$ and $T(d_1) = d$, then $d\geq a\vee b, \{a_1\}\vee\{b_1\}\subseteq T^{-1}(a\vee b)$, and $[a_1\vee b_1]\leq[d_1]$ (in view of the order relation introduced in first part) so that there are $a_2\in T^{-1}(a), b_2\in T^{-1}(b)$ such that $d_1\geq a_2\vee b_2\equiv a_1\vee b_1$(mod $\theta$).
Similarly, if $T^{-1}(a)\neq T^{-1}(b)$, $a_1 \in T^{-1}(a), b_1 \in T^{-1}(b)$, $d_1 \leq a_1$ and $d_1\leq b_1$ then there are $a_2\in T^{-1}(a), b_2\in T^{-1}(b)$, such that   $d_1\leq a_2\wedge b_2\equiv a_1\wedge b_1$ (mod $\theta$). This completes the proof of the theorem.

\end{document}